\documentclass[a10paper,english,10pt]{article}
\usepackage[T1]{fontenc}
\usepackage{babel}
\usepackage{graphicx,color}
\usepackage{babel}
\usepackage{amsfonts}
\usepackage{amssymb, amsmath, amsfonts, amsthm, mathrsfs}
\usepackage{hhline, latexsym}
\usepackage{float}
\usepackage[margin=2.5cm]{geometry}
\begin{document}
\author{Aymen Braghtha}
\title{Local monomialization conjecture of a singular foliation of Darboux type}
\maketitle
\begin{abstract}
After the nice result introduced by Belotto in [1] concerning the local monomialization of a singular foliation given by $n$ first integrals, this work is a continuation in the same spirit. In this paper, we introduce a important conjecture about local monomialization of a singular foliation of Darboux type (see section 1). This conjecture can be used to study pseudo-abelian integrals [2,4].
\end{abstract}
\section{Introduction}

Let $M$ be an analytic manifold of dimension $n+2$. Given a families of first integrals of Darboux type $H_{\epsilon}$ 
\begin{equation}
H_{\epsilon}(x,y)=H(x,y,\epsilon_1,\ldots,\epsilon_n)=\prod^{k}_{i=1}P^{a_i}(x,y,\epsilon_1,\ldots,\epsilon_n),\quad a_i>0.
\end{equation} 

Let $F$ be the foliation of codimension one in $M$ with coordinates $(x,y,\epsilon_1,\ldots,\epsilon_n)$ which is given by the analytic one form $\omega$
\begin{equation}
\omega=\frac{H_x}{\phi}dx+\frac{H_y}{\phi}dy+\sum^{n}_{i=1}\frac{H_{\epsilon_i}}{\phi}d\epsilon_i=0,
\end{equation}
where $H_x=\frac{\partial H}{\partial x}, H_y=\frac{\partial H}{\partial y}, H_{\epsilon_i}=\frac{\partial H}{\partial\epsilon_i}$ and $\phi=\prod^k_{i=1}P^{a_i-1}_i(x,y,\epsilon_1,\ldots,\epsilon_n)$ (integrating factor).

Let $F_i, i=1,\ldots,n$ are foliations of codimension one in $M$ with coordinates $(x,y,\epsilon_1,\ldots,\epsilon_n)$ which are given by the one forms $\omega_i$
 \begin{equation}
\omega_i=d\epsilon_i=0, \quad i=1,\ldots,n.
\end{equation}

Let $\mathcal{F}=(F,F_1,\ldots,F_n)$ be the result foliation of dimension one in $M$ where its leaves are given by the transversal intersection of leaves of $F, F_1,\ldots,F_n$. Otherwise speaking, the singular foliation $\mathcal{F}$ is given by 
\begin{align}
\Omega&=\omega\wedge\omega_1\wedge\ldots\wedge\omega_n\\
&=Q_1(x,y,\epsilon_1,\ldots,\epsilon_n)dx\wedge d\epsilon_1\wedge\ldots\wedge d\epsilon_n+Q_2(x,y,\epsilon_1,\ldots,\epsilon_n)dy\wedge d\epsilon_1\wedge\ldots\wedge d\epsilon_n=0,
\end{align}
where $Q_1=\frac{H_x}{\phi},Q_2=\frac{H_y}{\phi}$ are polynomials.

We shall say that $\Omega$ is a foliation of Darboux type with first integrals $(H,\epsilon_1,\ldots,\epsilon_n)$.\\
\linebreak
\textbf{Example.} Let $H_{\epsilon}(x,y)=H(x,y,\epsilon)=(x-\epsilon)^{a_1}(x-y)^{a_2}(x+y)^{a_3}$ be a the first integral of Darboux type. The foliation $F$ of codimension one in three dimensional space $M$ with coordinates $(x,y,\epsilon)$ is given by the one form
\begin{align*}
\omega&=(a_1(x-y)(x+y)+a_2(x-\epsilon)(x+y)+a_3(x-\epsilon)(x-y))dx\\
&-(a_2(x-\epsilon)(x+y)-a_3(x-\epsilon)(x-y))dy-a_1(x-y)(x+y)d\epsilon=0
\end{align*}
and the foliation $F_1$ of codimension one in $M$ is given by the one form
$$
\omega_1=d\epsilon=0.
$$
The result foliation $\mathcal{F}=(F,F_1)$ is given by the two-form 
$$
\Omega=\omega\wedge\omega_1=(a_1(x-y)(x+y)+a_2(x-\epsilon)(x+y)+a_3(x-\epsilon)(x-y))dx\wedge d\epsilon-(a_2(x-\epsilon)(x+y)-a_3(x-\epsilon)(x-y))dy\wedge d\epsilon=0
$$

Observe that the foliation $\mathcal{F}=(F,F_1)$ has a complicated singularity at the origin $(0,0,0)\subset D_0=\{\epsilon=0\}$.\\
\linebreak
\textbf{Conjecture.} \textsc{There exist sequences of local blowings-up such that the total transform of the foliation $\mathcal{F}: \omega\wedge\omega_1\wedge\ldots\wedge\omega_n=0$ has locally $n+1$ monomial first integrals $(z^{\gamma_0},z^{\gamma_1}...,z^{\gamma_n})$ where $z^{\gamma_i}=z_1^{\gamma_{i,1}}\cdots z_{n+2}^{\gamma_{i,n+2}}$ and the exponents matrix 
$$
m(a_1,\ldots,a_k)=
\begin{pmatrix}
\gamma_0^1&\ldots&\gamma_0^{n+2}\\
\gamma_1^1&\ldots&\gamma_1^{n+2}\\
\vdots&\vdots&\vdots\\
\gamma_n^1&\ldots&\gamma_n^{n+2}
\end{pmatrix}
$$
has a maximal rank.}
\section{Blowing-up of the foliation $\mathcal{F}$}
In this section, we introduce the fundamental idea to prove the conjecture which is based in first step on Hironaka's reduction of singularities [3].
Let $D_0=\{\epsilon_1=\epsilon_2=\ldots=\epsilon_n=0\}$ be a initial exceptional divisor.\\
\linebreak
\textbf{Theorem 1.} \emph{There exist a morphism $\Phi$ such that the pull-back foliation $\widetilde{\Phi}^{\ast}\mathcal{F}=\widetilde{\mathcal{F}}$ is given locally in neighborhood $U_1$ of the divisor $\widetilde{\Phi}^{\ast}(D_0)$ with coordinates $z=(z_1,\ldots,z_{n+2})$ by the following system
\begin{eqnarray}
\left\lbrace
\begin{array}{l}
\widetilde{H}=z^{\gamma_0}.\Delta_0,\\
\tilde{\epsilon}_1=z^{\gamma_1}.\Delta_1,\\
\vdots\\
\tilde{\epsilon}_n=z^{\gamma_n}.\Delta_n,
\end{array}
\right.
\end{eqnarray}
where $\Delta_i, i=0,\ldots,n$ are a units. }
\begin{proof}
(1) In first step, we monomialize the principal ideal $I_1=<P_1>$, Hironaka theorem's guarantee the existence oft a sequence of blow-ups $\widetilde{\Phi}_1=\widetilde{\Phi}^1_{n_1}\circ\widetilde{\Phi}^1_{n_1-1}\circ\ldots\circ\widetilde{\Phi}^1_{1}$ with initial center $C_0\subset D_0$ (which is possibly a submanifold of $M$) such that 
$$
(\widetilde{\Phi}^{\ast}_1P_1)^{a_1}=\delta_1\prod^{n+2}_{i=1}z^{a_1\tilde{\beta}_{i}^1}_{i},\quad \delta_1(0)\neq0.
$$
(2) In the second step, we consider the principal ideal $I_2=<\widetilde{\Phi}^{\ast}_1P_2>$ and by Hironaka theorem's there exist a sequence of blow-ups $\widetilde{\Phi}_2=\Phi^2_{n_2}\circ\Phi^2_{n_{2}-1}\circ\ldots\circ\Phi^2_{1}$ such that 
$$
(\widetilde{\Phi}^{\ast}_2\circ\widetilde{\Phi}^{\ast}_1P_2)^{a_2}=\delta_2\prod^{n+2}_{i=1}z^{a_2\tilde{\beta}_{i}^2}_{i},\quad\delta_2(0)\neq0.
$$
In the $k$-th step there exist a sequence of blow-ups $\widetilde{\Phi}_k=\Phi^k_{n_k}\circ\Phi^k_{n_{k}-1}\circ\ldots\circ\Phi^k_{1}$ such that the principal ideal $I_k=<\widetilde{\Phi}^{\ast}_{k-1}\circ\widetilde{\Phi}^{\ast}_{k-2}\circ\ldots\circ\widetilde{\Phi}^{\ast}_{1}P_k>$ has a normal crossings i.e. 
$$
(\widetilde{\Phi}^{\ast}_{k-1}\circ\widetilde{\Phi}^{\ast}_{k-2}\circ\ldots\circ\widetilde{\Phi}^{\ast}_{1}P_k)^{a_k}=\delta_k\prod^{n+2}_{i=1}z^{a_k\tilde{\beta}_{i}^k}_{i},\quad\delta_k(0)\neq0.
$$
Finally, after desingularisation of each polynomial $P_i$ of the first integral $H=\prod^{k}_{i=1}P_i^{a_i}$, the equations $z_1=0, \ldots, z_{n+2}=0$ are corresponding the irreducibles components of the exceptional divisor. For this raison after desingularisation of $\widetilde{\Phi}^{\ast}_{i-1}\circ\widetilde{\Phi}^{\ast}_{i-2}\circ\ldots\circ\widetilde{\Phi}^{\ast}_{1}P_i$, the polynomial $\widetilde{\Phi}_i^{\ast}\circ\widetilde{\Phi}^{\ast}_{i-1}\circ\widetilde{\Phi}^{\ast}_{i-2}\circ\ldots\circ\widetilde{\Phi}^{\ast}_{1}P_{i-1}$ has a normal crossings. So locally we have
\begin{eqnarray*}
\left\lbrace
\begin{array}{l}
\widetilde{H}=z^{\gamma_0}.\Delta_0,\\
\tilde{\epsilon}_1=z^{\gamma_1}.\Delta_1,\\
\vdots\\
\tilde{\epsilon}_n=z^{\gamma_n}.\Delta_n,
\end{array}
\right.
\end{eqnarray*}
where $z=(z_1,\ldots,z_{n+2}), \gamma_0=\sum^{k}_{i=1}a_i\beta_i, \beta_i=(\beta_i^1,\ldots,\beta_i^{n+2}), \gamma_i=(\gamma_i^1,\ldots,\gamma_i^{n+2})$.
\end{proof}
To complete the proof its necessairy to eliminate the units $\Delta_0, \Delta_1,\ldots,\Delta_n$ in the system (6). Now we define the resonant locus of the foliation $(z^{\gamma_0}.\Delta_0,z^{\gamma_1}.\Delta_1,\ldots,z^{\gamma_n}.\Delta_n)$
$$\mathcal{R}:=\{a=(a_1,\ldots,a_k): \gamma_0\wedge\sum^n_{j=1}\gamma_j=0\}.$$
To prove the conjecture, we distinguish two cases
\begin{itemize}
\item \textbf{generic case $a\notin\mathcal{R}$}.
\item \textbf{nongeneric case $a\in\mathcal{R}$}. 
\end{itemize}
\section{Some examples in dimension three}
To more understand the problem, we see some examples in dimension three.\\
\linebreak
\textbf{Example 1:} Let $\mathcal{F}$ be the local foliation which is obtained by after $k$ blow-ups. The foliation $\mathcal{F}$ is given by the following system
\begin{eqnarray} 
\left\lbrace
\begin{array}{l}
H_{a}=x^{a_1}y^{a_2}(1+z)\\
f=xy 
\end{array}
\right.
\end{eqnarray}
In this example we have  $\gamma_0: a_1\beta_1+a_2\beta_2$ where $\beta_1=(1,0,0), \beta_{2}=(0,1,0)$, $\gamma_1=(1,1,0)$ and $\mathcal{R}=\{a=(a_1,a_2): \gamma_0\wedge\gamma_1=0\}$. Our goal is to kill the unit $1+z$ in the first integral $H_a$ without modifying the second monomial $f$ in the sense to preserve its monomiality structure. For this raison, we distinguish two different cases:\\
\linebreak
(a) The generic case $a_1\neq -a_2\Leftrightarrow a\notin\mathcal{R}$: We take the change of variables $\tilde{x}=x(1+z)^{\frac{1}{a_1+a_2}}, \tilde{y}=y(1+z)^{\frac{1}{a_2+a_1}}$ and $\tilde{z}=z$. Then, we obtain the following system
\begin{eqnarray}
\left\lbrace
\begin{array}{l}
H_a=\tilde{x}^{a_1}\tilde{y}^{a_2}\\
f=\tilde{x}\tilde{y}
\end{array}
\right.
\end{eqnarray}\\
\linebreak
\textbf{Question:} How to calculate the generator vector field of the monomial foliation (7)?\\
\linebreak
Let us assume that the foliation $\mathcal{F}$ is generated locally by the vector field  $X(\tilde{x},\tilde{y},z)=\alpha_1\tilde{x}\frac{\partial}{\partial\tilde{x}}+\alpha_2\tilde{y}\frac{\partial}{\partial\tilde{y}}+\alpha_3z\frac{\partial}{\partial z}$ which satisfies   
$$
X(H)=X(\tilde{x}^{a_1}\tilde{y}^{a_2})=0,\qquad X(f)=X(\tilde{x}\tilde{y})=0
$$
To determine the vector $\alpha=(\alpha_1,\alpha_2,\alpha_3)$ we use the two following relations
$$
<\alpha,\gamma_0>=0\quad (\text{i.e.}\quad X(H)=X(\tilde{x}^{a_1}\tilde{y}^{a_2})=0),\quad<\alpha,\gamma_1>=0 \quad(\text{i.e.} \quad X(f)=X(\tilde{x}\tilde{y})=0).
$$
where $<,>$ is scalar product in $\mathbb{C}^3$. Finally, the vector $(\alpha_1,\alpha_2,\alpha_3)\in\{e_3\}$ and then
$$
\mathcal{F}=\{z\frac{\partial}{\partial z}\}
$$
Now, we express the vector field $X$ in the original coordinates $(x,y,z)$. If we write $X(x,y,z)=Ax\frac{\partial}{\partial x}+By\frac{\partial}{\partial y}+z\frac{\partial}{\partial z}$, to determine $A, B$ we use the fact that 
$$X(xy)=0\Longleftrightarrow A=-B$$
and so $X(x,y,z)=A(x\frac{\partial}{\partial x}-y\frac{\partial}{\partial y})+z\frac{\partial}{\partial z}$ on the other hand we have 
$$
Ax=X(x)=X(\tilde{x}(1+z)^{\frac{1}{a_1+a_2}})=z\frac{\partial}{\partial z}(\tilde{x}(1+z)^{\frac{1}{a_1+a_2}})=\tilde{x}(1+z)^{\frac{1}{a_1+a_2}-1}\frac{z}{a_1+a_2}
$$
Finally, we obtain 
$$
X(x,y,z)=\frac{1}{a_1+a_2}\frac{z}{1+z}(x\frac{\partial}{\partial x}-y\frac{\partial}{\partial y})+z\frac{\partial}{\partial z}\Rightarrow Y(x,y,z)=x\frac{\partial}{\partial x}-y\frac{\partial}{\partial y}+(a_1+a_2)(1+z)\frac{\partial}{\partial z}
$$
\textbf{Remark 1.} \emph{In dimension three, if we consider the foliation $\mathcal{F}$ which is given locally by 
$$
\left\lbrace
\begin{array}{l}
f_1=x^{a}y^{b}z^{c}\\
f_2=x^{\tilde{a}}y^{\tilde{b}}z^{\tilde{c}}
\end{array}
\right.
$$
where rank$\begin{pmatrix}
a&b&c\\
\tilde{a}&\tilde{b}&\tilde{c}
\end{pmatrix}$=2. The generator vector field $X$ of the form 
$$X(x,y,z)=\hat{a}x\frac{\partial}{\partial x}+\hat{b}y\frac{\partial}{\partial y}+\hat{c}z\frac{\partial}{\partial z},$$
where 
$$<(\hat{a},\hat{b},\hat{c}),(a,b,c)>=0,\quad\text{and}\quad<(\hat{a},\hat{b},\hat{c}),(\tilde{a},\tilde{b},\tilde{c})>=0.$$
}\\
\linebreak  
In our example we observe that in the neighborhood of the leaf $\{z=0\}$ the vector field 
$$
Y\simeq x\frac{\partial}{\partial x}-y\frac{\partial}{\partial y}+(a_1+a_2)z\frac{\partial}{\partial z}
$$  
is linearizable and consequently $Y$ is transversal to the leaf $\{z=0\}$.\\
(b) The problem suppose where $a_1=-a_2$ i.e $a\in\mathcal{R}=\{a=(a_1,a_2): \gamma_0\wedge\gamma_1=0\}$. In this case near the leaf $\{z=0\}$, we have 
$$Y\simeq x\frac{\partial}{\partial x}-y\frac{\partial}{\partial y}.$$
It is clear that the condition of transversality of $Y$ and the leaf $\{z=0\}$ is not satisfieted.\\
\linebreak
\textbf{Example 2:} let $\mathcal{F}$ be the local foliation which is given (after a sequence of blow-ups) by 
$$
\left\lbrace
\begin{array}{l}
H=x^{a_1}y^{a_2}z^{a_3}(1+g(x,y,z))\\
f=xyz.
\end{array}
\right.
$$
The foliation $\mathcal{F}$ is given also 
$$
\left\lbrace
\begin{array}{l}
\frac{H}{f^{a_1}}=y^{a_2-a_1}z^{a_3-a_1}(1+g(x,y,z))\\
f=xyz.
\end{array}
\right.
$$
If $a=(a_1,a_2,a_3)\notin\mathcal{R}=\{a: (a_1\beta_1+a_2\beta_2+a_3\beta_3)\wedge (1,1,1)=0\}$ (resonant locus), we can take the following variables change $x=\tilde{x}, y=\tilde{y}(1+g(x,y,z))^{\frac{1}{a_2-a_3}}$ and $z=\tilde{z}(1+g(x,y,z))^{\frac{1}{a_3-a_2}}$
and in this case the local foliation $\mathcal{F}$ is generated by the vector field 
$$
X(\tilde{x},\tilde{y},\tilde{z})=(a_2-a_3)\tilde{x}\frac{\partial}{\partial \tilde{x}}+(a_3-a_1)\tilde{y}\frac{\partial}{\partial\tilde{y}}+(a_1-a_2)\tilde{z}\frac{\partial}{\partial\tilde{z}}
$$
let us express the vector field $X$ in the original coordinates $(x,y,z)$, so we have
$$
X(x)=X(\tilde{x})=(a_2-a_3)x\frac{\partial}{\partial x}
$$
 
$$
X(\tilde{y}(1+g(x,y,z)^{\frac{1}{a_2-a_3}}))=(a_3-a_1)y+\frac{1}{a_2-a_3}y \frac{X(g(x,y,z))}{1+g(x,y,z)}
$$
and
$$
X(\tilde{z}(1+g(x,y,z)^{\frac{1}{a_3-a_2}}))=(a_1-a_2)z+\frac{1}{a_3-a_2}z \frac{X(g(x,y,z))}{1+g(x,y,z)}.
$$
Finally the vector field $X$ of the form 
$$
X(x,y,z)=(a_2-a_3)x\frac{\partial}{\partial x}+(a_3-a_1)y\frac{\partial}{\partial y}+(a_1-a_2)z\frac{\partial}{\partial z}+\frac{1}{a_2-a_3} \frac{X(g(x,y,z))}{1+g(x,y,z)}(y\frac{\partial}{\partial x}-z\frac{\partial}{\partial z}).
$$
\textbf{Proposition 1.} \emph{For $a\notin\mathcal{R}$, there exist a local diffeomorphism $\phi: z=(z_1,\cdots,z_{n+2})\mapsto\tilde{z}=(\tilde{z}_1,\cdots,\tilde{z}_{n+2})$ such that the foliation $\mathcal{F}$ is given locally by $( \tilde{z}^{\gamma_0},\cdots,\tilde{z}^{\gamma_{n+2}})$}
\begin{proof}
We just make a convenable change of variables.
\end{proof}

\textbf{Open question.} To complete the proof of Conjecture we must solve the nongeneric case $a=(a_1,\ldots,a_k)\in\mathcal{R}$ because in this case the rank of exponent matrix 
$$
m(a_1,\ldots,a_k)=
\begin{pmatrix}
\gamma_0^1&\ldots&\gamma_0^{n+2}\\
\gamma_1^1&\ldots&\gamma_1^{n+2}\\
\vdots&\vdots&\vdots\\
\gamma_n^1&\ldots&\gamma_n^{n+2}
\end{pmatrix}
$$
is not maximal.

\textsl{Universit\'{e} de Bourgogne, Institut de Math\'{e}matiques de\\
 Bourgogne, U.M.R. 5584 du C.N.R.S., B.P. 47870, 21078 Dijon \\
cedex - France.\\}
E-mail adress: aymenbraghtha@yahoo.fr

\end{document}